\documentclass{article}

\usepackage{arxiv}

\usepackage[utf8]{inputenc} % allow utf-8 input
\usepackage[T1]{fontenc}    % use 8-bit T1 fonts
\usepackage{hyperref}       % hyperlinks
\usepackage{url}            % simple URL typesetting
\usepackage{booktabs}       % professional-quality tables
\usepackage{amsfonts}       % blackboard math symbols
\usepackage{nicefrac}       % compact symbols for 1/2, etc.
\usepackage{microtype}      % microtypography
\usepackage{lipsum}
\usepackage{graphicx}
\usepackage{amsmath}
\usepackage{amsfonts}
\usepackage[ruled,vlined]{algorithm2e}
\usepackage{multirow}

\title{Restricted Boltzmann Machine Assignment Algorithm: \\ Application to solve many-to-one matching problems on weighted bipartite graph}

\author{
  Francesco Curia \\
  Department of Statistical Science\\
  Sapienza, University of Rome\\
  Rome, 00185 Italy \\
  \texttt{francesco.curia@uniroma1.it} \\
}

\begin{document}
\maketitle

\begin{abstract}
In this work an iterative algorithm based on unsupervised learning is presented, 
specifically on a Restricted Boltzmann Machine (RBM) to solve a perfect matching 
problem on a bipartite weighted graph. Iteratively 
is calculated the weights $ w_{ij} $ and the bias 
parameters $\theta = ( a_i, b_j) $ that maximize the energy function and assignment element $i$ to element $j$.
An application of real problem is presented to show the potentiality of this algorithm.
\end{abstract}


\section{Introduction}
\label{intro}

The assignment problems fall within the combinatorial optimization problems,
the problem of matching on a bipartite weighted graph is one of the major problems faced
in this compound. Numerous resolution methods and algorithms have been proposed in recent times and many
have provided important results, among them for example we find: Constructive heuristics, Meta-heuristics,
Approximation algorithms, Iper-heuristics, and other methods. 
Combinatorial optimization deals with finding the 
optimal solution between the collection of finite possibilities. 
The finished set of possible solutions. The heart of the problem of 
finding solutions in combinatorial optimization is based on the efficient 
algorithms that present with a polynomial computation time in the input dimension. 
Therefore, when dealing with certain combinatorial optimization problems one must ask
with what speed it is possible to find the solutions or the optimal problem 
solution and if it is not possible to find a resolution method of this type, 
which approximate methods can be used in polynomial computational times that lead to stable explanations.
Solve this kind of problem in polynomial time $o(n)$
has long been the focus of research in this area until Edmonds [1] 
developed one of the most efficient methods. 
Over time other algorithms have been developed, for example the fastest of them is the Micali e Vazirani algorithm
[2], Blum [3] and Gabow and Tarjan [4]. The first
of these methods is an improvement on that of Edmonds, 
the other algorithms use different logics, but all of them
with computational time equal to $o(m \sqrt n)$.
The problem is fundamentally the following: we imagine a situation in which respect for the characteristics detected on a given phenomenon is to be assigned between elements of two sets, as for example in one of the most known problems such as the tasks and the workers to be assigned to them. A classical maximum cardinality matching algorithm to take the maximum weight range and assign it, in a decision support system, through the domain expert this could also be acceptable, but in a totally automatic system like a system could be of artificial intelligence that puts together pairs of elements on the basis of some characteristics, this way would not be very reliable, totally removing the user's control. Another problem related to this kind of situation is that of features. Let's take as an example a classic problem of flight-gate assignment in an airport, on the basis of the history we could have information about the flight, the gates and the time, the flight number and maybe the airline. Little information that even through the best of feature enginering would lead to a model of machine learning, specifically of classification, very poor in information. Treating the same problem with classical optimization, as done so far, would lead to solving it with a perfect matching of maximum weight, and we would return to the beginning.

\section{Matching problems}
\label{sec:1}
Matching problems are among the fundamental problems in combinatorial optimization.
In this set of notes, we focus on the case when the underlying graph is bipartite.
We start by introducing some basic graph terminology. A graph $G = (V, E)$ consists of
a set $V = A \cup B$ of vertices and a set $E$ of pairs of vertices called edges. For an edge $e = (u, v)$, we
say that the endpoints of e are $u$ and $v$; we also say that $e$ is incident to $u$ and $v$. A graph
$G = (V, E)$ is bipartite if the vertex set $V$ can be partitioned into two sets $A$ and $B$ (the
bipartition) such that no edge in $E$ has both endpoints in the same set of the bipartition. A
matching M is a collection of edges such that every vertex of $V$ is incident to at most
one edge of $M$. If a vertex $v$ has no edge of $M$ incident to it then $v$ is said to be exposed
(or unmatched). A matching is perfect if no vertex is exposed; in other words, a matching is
perfect if its cardinality is equal to $|A| = |B|$. In the literature several examples of the 
real world have been treated such as the assignment of children to certain schools [5], 
or as donors and patients [6] and workers at companies [7].
The problem of the weighted bipartite matching finds the feasible match with the maximum available weight.
This problem was developed in several areas, such as in the work of [8]
about protein and structure alignment, or within the computer
vision as documented in the work of [9] or as in the paper by [10] 
in which the similarity of the texts is estimated. Other jobs have faced
this problem in the classification [11],[12] e [13], but not for many to one correspondence.
The mathematical formulation can be solved by presenting it as a linear program. 
Each edge $(i,j)$, where $i$ is in $A$ and $j$ is in $B$, has a weight $w_{ij}$. 
For each edge $(i,j)$  we have a decision variable

\begin{equation}
 x_{ij} =\begin{cases} 1 & \mbox{if the edge is contained in the matching} \\ 0 & \mbox{otherwise}
\end{cases} 
\end{equation}

and  $x_{ij}\in \mathbb {Z} {\text{ for }}i,j \in A,B$, and we have the following LP:

\begin{equation}
\begin{aligned}
& \underset{x_{ij}}{\text{max}} \sum _{(i,j)\in A\times B}w_{ij}x_{ij}
\end{aligned}
\end{equation}

\begin{equation}
\begin{aligned}
\sum _{j\in B}x_{ij}=1{\text{ for }}i\in A
\end{aligned}
\end{equation}

\begin{equation}
\begin{aligned}
\sum _{i\in A}x_{ij}=1{\text{ for }}j\in B
\end{aligned}
\end{equation}

\begin{equation}
\begin{aligned}
0 \leq x_{ij} \leq 1{\text{ for }}i,j \in A,B
\end{aligned}
\end{equation}

\begin{equation}
\begin{aligned}
x_{ij}\in \mathbb {Z} {\text{ for }}i,j\in A,B
\end{aligned}
\end{equation}

\begin{figure}
  \includegraphics[width=0.55\textwidth]{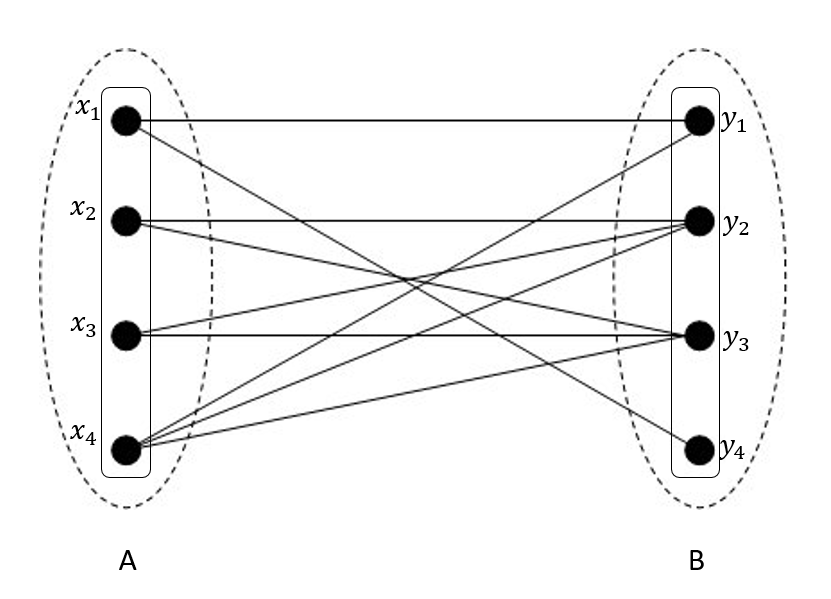}
\caption{Bipartite Weighted Matching}
\label{fig:1}       % Give a unique label
\end{figure}

%\cite{RefB} and \cite{RefJ}.
\section{Motivation}
\label{sec:2}

The problem in a weighted bipartite graph $G = (V = A \cup B, W)$, 
is when we have different weights (historical data) $W = \left(w_{11}, w_{12} ..., w_{ij} \right)$ for the set of nodes 
to which corresponds the same set of nodes $B$. One of most popular solution is Hungarian algorithm [16] . 
The assignment rule therefore  in many real cases could be misleading and limiting as well 
as it could be unrealistic as a solution. Machine learning (ML) algorithms are 
increasingly gaining ground in applied sciences such as engineering, biology, medicine, etc.
both supervised and unsupervised learning models. The matching problem 
in this case can be seen as a set of inputs $ x_1, ..., x_k$ (in our case the nodes of the set A) 
and a set of ouput $y_1, ..., y_k$ (the respective nodes of the set $B$), weighed by 
a series of weights $w_{11}, ... w_{ij}$, which inevitably recall the structure of a classic neural network. 
The problem is that in this case there would be a number of classes (in the case of assignment) 
equal to the number of inputs. Considering it as a classic machine learning problem, the difficulty would lie 
in the features and their engineering, on the one hand, while on the other 
the number of classes to predict (assign) would be very large.
For example, if we think about matching applicants and jobs, if we only had the name of a 
candidate for the job, we would have very little info to build a robust machine learning model, and even a 
good features engineering would not lead to much, but having available other information on the 
candidate they could be extracted it is used case never as "weight" to build a neural network, 
but even in this case the constraint of a classic optimization model solved with ML techniques would 
not be maintained, let's say we would "force" a little hand. While what we want to present in 
this work is the resolution of a classical problem of matching (assignment) through the application 
of a ML model, in this case of a neural network, which as already said maintains the mathematical 
structure of a node (input) and arc (weight) but instead of considering the output of arrival 
(the set $B$) as classification label (assignment) in this case we consider an unsupervised 
neural network, specifically an Restricted Boltzmann Machine.

\section{Contributions}
\label{sec:3}

The contributions of this work are mainly of two types: the first is related to the ability to use an unsupervised machine learning model to solve a classical optimization problem which in turn has the mathematical structure of a neural network based on two layers, in our case of RBM a visible and a hidden one. In this case the nodes of the set $ B $ become the variables of the model and the number of times that the node $i$ has been assigned to the node $j$ (for example in problems that concerns historical data analysis), becomes the weight $ w_ {ij} $ which at its turn becomes the value of the variable $i$-th in the RBM model. The second is the ability to solve real problems, as we will see later in the article, in which it is necessary to carry out a matching between elements of two sets and the maximum weight span is not said to be the best assignment, especially in the case where the problem is many-to-one, like many real problems.

\section{Restricted Boltzmann Machine}
\label{sec:3}

Restricted Boltzmann Machine is an unsupervised method neural networks based [14]; the algorithm learns one layer of hidden features. When the number of hidden units is smaller than that of visual units, the
hidden layer can deal with nonlinear complex dependency and structure of data, capture deep relationship from input data , and represent the input data more compactly.
Assuming there are $c$ visible units and $m$ hidden units in an Restricted Boltzmann Machine. So $v_i$ for $i = 1,...,c$ indicates the state of the $i-$th visible
unit, where 

\begin{equation}
 v_i =\begin{cases} 1 & \mbox{if the i-th term is annotated to the element } \\ 0 & \mbox{otherwise}
\end{cases} 
\end{equation}

for $i=1,...,c$ and furthemore we have 

\begin{equation}
h_j =\begin{cases} 1 & \mbox{the state of hidden unit is active} \\ 0 & \mbox{otherwise}
\end{cases} 
\end{equation}

for  $j=1,...,m$ and $w_{ij}$ is  the weight associated with the connection between $v_i$ and $h_j$ and define also the joint configuration $(v,h)$.

\begin{figure*}[ht]
  \includegraphics[width= 0.75 \textwidth]{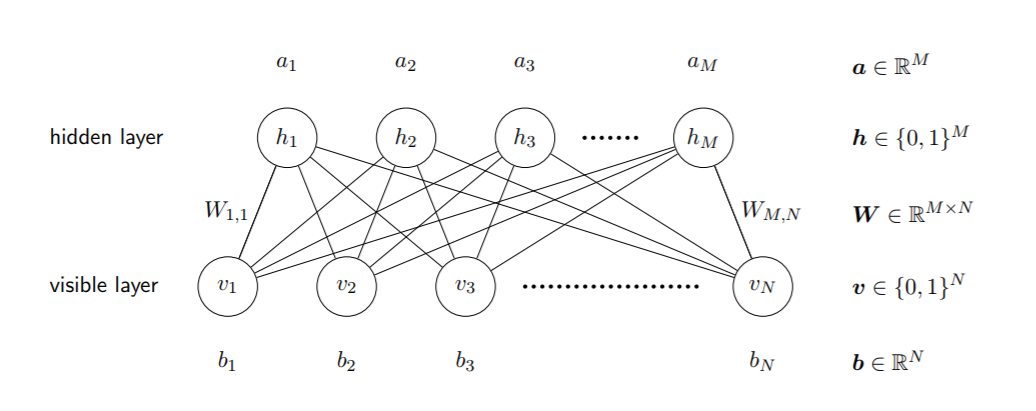}
\caption{Restricted Boltzamann Machine Structures}
\label{fig:2}       % Give a unique label
\end{figure*}

The energy function that capturing the interaction patterns between visible layer and hidden layer is define as follow:
\begin{equation}
E(v,h| \theta ) = - \sum_{i = 1}^{c} a_i v_i - \sum_{j = 1}^{m} b_j h _j - \sum_{i = 1}^{c}  \sum_{j = 1}^{m} v_i h_j w_{ij}
\end{equation}

where $\theta = \left( w_{ij}, a_i, b_j\right)$ are parameters of the model: $a_i$ and $b_j$ are biases for the visible and hidden variables, respectively. The parameters $w_{ij}$ is  the weights of connection between visible variables and hidden variables. The joint probability is represented by the follow quantity:
\begin{equation}
p(v,h) = \frac{e^-E(v,h)}{Z}
\end{equation}

where

\begin{equation*}
Z = \sum_{v,h} e^{-E(v,h)}
\end{equation*}

is a normalization constant and the conditional distributions over the visible and hidden units are given
by sigmoid functions as follows:

\begin{equation}
p(v_i = 1 | h) = \sigma \left( \sum_{j = 1}^{m} w_{ij} h_j + a_i \right)
\end{equation}

\begin{equation}
p(h_j = 1 | v) = \sigma \left( \sum_{i = 1}^{c} w_{ij} v_i + b_j  \right)
\end{equation}

where $\sigma = \frac{1}{1 + e^{-x}}$

$\\$
RBM are trained to optimizie the product of probabilities assigned to some training set $V$ (a matrix, each row of which is treated as a visible vector  v)

\begin{equation}
\begin{aligned}
& \underset{w_{ij}}{\text{arg} \ \text{max}} \prod_{i=1}^c p(v_i)
\end{aligned}
\end{equation}

The RBM training takes place through the Contrastive Divergence Algorithm (see Hinton 2002 [15]).  For the (4) and (7) we can pass to log-likelihood formulation

\begin{equation}
L_v = log \sum_{h} e^{-E(v,h)} - log \sum_ {v,h} e^{-E(v,h)}
\end{equation}

and derivate the quantity 

\begin{equation}
\frac{\partial L_v}{\partial w_{ij}} = \sum_{h} p(h|v) \cdot v_i h_j - \sum_{v,h} p(v,h) \cdot v_i h_j
\end{equation}

\begin{equation}
\frac{\partial L_v}{\partial w_{ij}} = \mathbb{E} [p(h|v)]  - \mathbb{E}[p(v,h)]
\end{equation}
In the above expression, the first quantity represents the expectation of $v_i \cdot h_j$   to the conditional probability of the hidden states given the visible states and the second term represents the expectation of $v_i \cdot h_j$  to the joint probability of the visible and hidden states. In order to maximize the (7), which involves log of summation, there is no analytical solution and will use the stochastic gradient ascent technique. 
In order to compute the unknown values of the weights $w_{ij}$ such that they maximize the above likelihood, we compute using gradient ascent :
\begin{equation*}
w_{ij}^{k+1}  = w_{ij}^{k} + \alpha \cdot \frac{\partial L_v}{\partial w_{ij}} %- \beta \cdot w_{ij}^{k} + \lambda \cdot  \left(w_{ij}^{k}  - w_{ij}^{k-1}\right)
\end{equation*}
where $\alpha \in (0,1) $  is the learning rate. At this formulation we can add an term of penality and obtain the follow
\begin{equation}
w_{ij}^{k+1}  = w_{ij}^{k} + \alpha \cdot \frac{\partial L_v}{\partial w_{ij}} + \eta \cdot \frac{w_{ij}^{k-1}}{w_{ij}^k}
\end{equation}
where $\eta > 0$ is penalty parameters and $ \eta \cdot \frac{w_{ij}^{k-1}}{w_{ij}^k}$ measure the contribute of weights variation in the $k$-th step of update.
 %, the term $\beta \cdot w_{ij}^{k}$ is for regularization purpose so that the magnitude of the weights do not become too large and $\lambda \cdot \left(w_{ij}^{k}  - w_{ij}^{k-1}\right)$ is the momentum term to prevent   oscillations.

\section{Restricted Boltzmann Machine Assignment Algorithm (RBMAA)}
\label{sec:4}

In this section is presented the algorithm RBM based and is explaneid the single steps. 

\begin{enumerate}
\item In this step the algorithm takes as input the matrix $W$ relative to the weights of the assignments, either the weight $ w_ {ij} $ represents the number of times that the element $ i $ has been assigned to the element $ j $
\newline
\item The matrix $W$ is binarized. For each element $ w_ {ij} > $ 0  take 1 or 0 and an a matrix $\tilde W$  is created with elements 0-1.
\newline
\item In this step the RBM is applied taking as input the binary matrix $\tilde W $. The probability product related to the visible units of the RBM is maximized (13). The weights $ w_ {ij} $ are updated according to (17) and the biases according to the RBM training rule. Once these updated values are obtained, the optimized value $ \hat p (v_i) $ if this is greater than $\epsilon >$ 0 assigns 1 otherwise 0 and the iteration of the RBM is restarted until for each row of the matrix there is not a single value 1 which is the same as assignment. The quantile of level $ \alpha = $ 0.99 was used to determine the threshold $\epsilon$.
\newline
\item The output of the algorithm is a matrix with values 0-1 and on each line we obtain a single value equal to 1 which is equivalent to the assignment of the element $ i $ to the element $ j $.
\end{enumerate}

The pseudocode is presented in the next section.

\DontPrintSemicolon
\newcommand{\To}{\mbox{\upshape\bfseries to}}

\begin{algorithm}[t]
  \caption{Restricted Boltzmann Assignment}
  %\textbf{Phase 1}\;
  \nl \textbf{Input}:  \\
$\epsilon$, threshold value \\
 $b_j$, hidden units bias value\\
 $a_i$, visible units bias value \\
 $\alpha$, learning rate for visible bias updating \\
$\beta$, learning rate for hidden bias updating\\

 Matrix $n  \times m$, $W = \{ w_{ij} \}$ number of times which $i$ has been assigned to $j$\ 
\\ 
\quad

  \nl {% 
Binarize $W$ to get $\dot W = \{ \dot w_{ij} \} \in [0,1] $
  \;
    \ForEach{$\text{i} \ in \ W$}{
    \eIf{$w_{ij} >0 $ }{
   1\;
   \;
   }{
   0\;
  }
    }
 \;
    
 \nl {%
    \While{$\sum_{i = 1}^{c} \hat w_{ij} \ne $ 1}{\quad

$\underset{\dot w_{ij}}{\text{arg} \ \text{max}} \prod_{i=1}^c p(v_i)$
\\
{\textbf{Update} $w_{ij}$,
$w_{ij}^{k+1}  = w_{ij}^{k} + \alpha \cdot \frac{\partial L_v}{\partial w_{ij}}$

\textbf{Update}  $b_{j}$,
$b_{j}^{k+1}  = b_{j}^{k} + \beta \cdot \frac{\partial L_v}{\partial b_{j}} $

\textbf{Update}  $a_{i}$,
$a_{i}^{k+1}  = a_{i}^{k} + \alpha \cdot \frac{\partial L_v}{\partial a_{i}} $
}

\eIf{$\hat p(v_i) > \epsilon $}{1 \;
   \;
   }{
   0\;
  }
 }
 \;
}

\nl Output:  Matrix $n  \times m$, $\hat P = \{ 0,1 \}$ such that $\sum_{i} \hat p_{ij} = 1$, $\forall i \in A$

}

\end{algorithm}

\section{Application of the RBMAA  to a real problem}
\label{sec:5}
Now we proceed to provide the results of the application of the algorithm (see Appendix). The problem instance is 351 elements for the set $ A $ and 35 for the set $ B $; the goal is to assign for each element of $ A $, one and only one element of $ B $, so as to have the sum for row equal to 1 as in (3) and in step 4 of the algorithm. The weights $ w_ {ij} $ are represented by the number of times the element $ a_i \in A $ has been assigned to the element $ b_j \in B $, based on a set of historical data that we have relative to flight-gate assignments of a well-known international airport. The difficulty is the one discussed in the first part of the work, in which we want to obtain a robust machine learning algorithm that classifies and assigns the respective gate for each flight. Starting from the available features, the algorithm presented in the work was implemented. The computational results are very interesting in terms of calculation speed and assignment.

\section{Conclusions and discussion}
\label{sec:6}
This can be the starting point for more precise, fast and sophisticated algorithms, which combine the combinatorial optimization with machine learning but on the basis of unsupervised learning and not only on the optimization of cost functions.

$\\$
$\\$
$\\$
$\\$

$\\$
$\\$
$\\$
$\\$
$\\$
$\\$
$\\$
$\\$
$\\$
$\\$
$\\$
$\\$

\section{Appendix}

\begin{table}[htpb]
\resizebox{\textwidth}{!}{%
\begin{tabular}{llllllllllll|l|l|}
\hline
\multicolumn{1}{|l|}{\textbf{Node1}} & \multicolumn{1}{l|}{\textbf{Node2}} & \multicolumn{1}{l|}{\textbf{Node1}} & \multicolumn{1}{l|}{\textbf{Node2}} & \multicolumn{1}{l|}{\textbf{Node1}} & \multicolumn{1}{l|}{\textbf{Node2}} & \multicolumn{1}{l|}{\textbf{Node1}} & \multicolumn{1}{l|}{\textbf{Node2}} & \multicolumn{1}{l|}{\textbf{Node1}} & \multicolumn{1}{l|}{\textbf{Node2}} & \multicolumn{1}{l|}{\textbf{Node1}} & \textbf{Node2} & \textbf{Node1} & \textbf{Node2} \\ \hline
\multicolumn{1}{|l|}{A2}             & \multicolumn{1}{l|}{B1}             & \multicolumn{1}{l|}{A173}           & \multicolumn{1}{l|}{B5}             & \multicolumn{1}{l|}{A203}           & \multicolumn{1}{l|}{B10}            & \multicolumn{1}{l|}{A268}           & \multicolumn{1}{l|}{B15}            & \multicolumn{1}{l|}{A158}           & \multicolumn{1}{l|}{B20}            & \multicolumn{1}{l|}{A48}            & B25            & A183           & B30            \\ \hline
\multicolumn{1}{|l|}{A72}            & \multicolumn{1}{l|}{B1}             & \multicolumn{1}{l|}{A208}           & \multicolumn{1}{l|}{B5}             & \multicolumn{1}{l|}{A238}           & \multicolumn{1}{l|}{B10}            & \multicolumn{1}{l|}{A303}           & \multicolumn{1}{l|}{B15}            & \multicolumn{1}{l|}{A193}           & \multicolumn{1}{l|}{B20}            & \multicolumn{1}{l|}{A83}            & B25            & A218           & B30            \\ \hline
\multicolumn{1}{|l|}{A107}           & \multicolumn{1}{l|}{B1}             & \multicolumn{1}{l|}{A243}           & \multicolumn{1}{l|}{B5}             & \multicolumn{1}{l|}{A308}           & \multicolumn{1}{l|}{B10}            & \multicolumn{1}{l|}{A338}           & \multicolumn{1}{l|}{B15}            & \multicolumn{1}{l|}{A228}           & \multicolumn{1}{l|}{B20}            & \multicolumn{1}{l|}{A93}            & B25            & A253           & B30            \\ \hline
\multicolumn{1}{|l|}{A117}           & \multicolumn{1}{l|}{B1}             & \multicolumn{1}{l|}{A278}           & \multicolumn{1}{l|}{B5}             & \multicolumn{1}{l|}{A318}           & \multicolumn{1}{l|}{B10}            & \multicolumn{1}{l|}{A22}            & \multicolumn{1}{l|}{B16}            & \multicolumn{1}{l|}{A263}           & \multicolumn{1}{l|}{B20}            & \multicolumn{1}{l|}{A118}           & B25            & A288           & B30            \\ \hline
\multicolumn{1}{|l|}{A142}           & \multicolumn{1}{l|}{B1}             & \multicolumn{1}{l|}{A313}           & \multicolumn{1}{l|}{B5}             & \multicolumn{1}{l|}{A343}           & \multicolumn{1}{l|}{B10}            & \multicolumn{1}{l|}{A57}            & \multicolumn{1}{l|}{B16}            & \multicolumn{1}{l|}{A273}           & \multicolumn{1}{l|}{B20}            & \multicolumn{1}{l|}{A153}           & B25            & A323           & B30            \\ \hline
\multicolumn{1}{|l|}{A177}           & \multicolumn{1}{l|}{B1}             & \multicolumn{1}{l|}{A348}           & \multicolumn{1}{l|}{B5}             & \multicolumn{1}{l|}{A27}            & \multicolumn{1}{l|}{B11}            & \multicolumn{1}{l|}{A92}            & \multicolumn{1}{l|}{B16}            & \multicolumn{1}{l|}{A298}           & \multicolumn{1}{l|}{B20}            & \multicolumn{1}{l|}{A188}           & B25            & A7             & B31            \\ \hline
\multicolumn{1}{|l|}{A212}           & \multicolumn{1}{l|}{B1}             & \multicolumn{1}{l|}{A32}            & \multicolumn{1}{l|}{B6}             & \multicolumn{1}{l|}{A75}            & \multicolumn{1}{l|}{B11}            & \multicolumn{1}{l|}{A127}           & \multicolumn{1}{l|}{B16}            & \multicolumn{1}{l|}{A333}           & \multicolumn{1}{l|}{B20}            & \multicolumn{1}{l|}{A223}           & B25            & A42            & B31            \\ \hline
\multicolumn{1}{|l|}{A247}           & \multicolumn{1}{l|}{B1}             & \multicolumn{1}{l|}{A67}            & \multicolumn{1}{l|}{B6}             & \multicolumn{1}{l|}{A97}            & \multicolumn{1}{l|}{B11}            & \multicolumn{1}{l|}{A162}           & \multicolumn{1}{l|}{B16}            & \multicolumn{1}{l|}{A17}            & \multicolumn{1}{l|}{B21}            & \multicolumn{1}{l|}{A258}           & B25            & A77            & B31            \\ \hline
\multicolumn{1}{|l|}{A282}           & \multicolumn{1}{l|}{B1}             & \multicolumn{1}{l|}{A102}           & \multicolumn{1}{l|}{B6}             & \multicolumn{1}{l|}{A110}           & \multicolumn{1}{l|}{B11}            & \multicolumn{1}{l|}{A197}           & \multicolumn{1}{l|}{B16}            & \multicolumn{1}{l|}{A52}            & \multicolumn{1}{l|}{B21}            & \multicolumn{1}{l|}{A293}           & B25            & A147           & B31            \\ \hline
\multicolumn{1}{|l|}{A292}           & \multicolumn{1}{l|}{B1}             & \multicolumn{1}{l|}{A137}           & \multicolumn{1}{l|}{B6}             & \multicolumn{1}{l|}{A132}           & \multicolumn{1}{l|}{B11}            & \multicolumn{1}{l|}{A232}           & \multicolumn{1}{l|}{B16}            & \multicolumn{1}{l|}{A87}            & \multicolumn{1}{l|}{B21}            & \multicolumn{1}{l|}{A328}           & B25            & A182           & B31            \\ \hline
\multicolumn{1}{|l|}{A317}           & \multicolumn{1}{l|}{B1}             & \multicolumn{1}{l|}{A172}           & \multicolumn{1}{l|}{B6}             & \multicolumn{1}{l|}{A202}           & \multicolumn{1}{l|}{B11}            & \multicolumn{1}{l|}{A267}           & \multicolumn{1}{l|}{B16}            & \multicolumn{1}{l|}{A122}           & \multicolumn{1}{l|}{B21}            & \multicolumn{1}{l|}{A12}            & B26            & A217           & B31            \\ \hline
\multicolumn{1}{|l|}{A1}             & \multicolumn{1}{l|}{B2}             & \multicolumn{1}{l|}{A207}           & \multicolumn{1}{l|}{B6}             & \multicolumn{1}{l|}{A237}           & \multicolumn{1}{l|}{B11}            & \multicolumn{1}{l|}{A302}           & \multicolumn{1}{l|}{B16}            & \multicolumn{1}{l|}{A157}           & \multicolumn{1}{l|}{B21}            & \multicolumn{1}{l|}{A47}            & B26            & A252           & B31            \\ \hline
\multicolumn{1}{|l|}{A11}            & \multicolumn{1}{l|}{B2}             & \multicolumn{1}{l|}{A242}           & \multicolumn{1}{l|}{B6}             & \multicolumn{1}{l|}{A307}           & \multicolumn{1}{l|}{B11}            & \multicolumn{1}{l|}{A337}           & \multicolumn{1}{l|}{B16}            & \multicolumn{1}{l|}{A167}           & \multicolumn{1}{l|}{B21}            & \multicolumn{1}{l|}{A82}            & B26            & A287           & B31            \\ \hline
\multicolumn{1}{|l|}{A36}            & \multicolumn{1}{l|}{B2}             & \multicolumn{1}{l|}{A277}           & \multicolumn{1}{l|}{B6}             & \multicolumn{1}{l|}{A342}           & \multicolumn{1}{l|}{B11}            & \multicolumn{1}{l|}{A21}            & \multicolumn{1}{l|}{B17}            & \multicolumn{1}{l|}{A192}           & \multicolumn{1}{l|}{B21}            & \multicolumn{1}{l|}{A152}           & B26            & A322           & B31            \\ \hline
\multicolumn{1}{|l|}{A71}            & \multicolumn{1}{l|}{B2}             & \multicolumn{1}{l|}{A312}           & \multicolumn{1}{l|}{B6}             & \multicolumn{1}{l|}{A26}            & \multicolumn{1}{l|}{B12}            & \multicolumn{1}{l|}{A56}            & \multicolumn{1}{l|}{B17}            & \multicolumn{1}{l|}{A227}           & \multicolumn{1}{l|}{B21}            & \multicolumn{1}{l|}{A187}           & B26            & A6             & B32            \\ \hline
\multicolumn{1}{|l|}{A106}           & \multicolumn{1}{l|}{B2}             & \multicolumn{1}{l|}{A347}           & \multicolumn{1}{l|}{B6}             & \multicolumn{1}{l|}{A61}            & \multicolumn{1}{l|}{B12}            & \multicolumn{1}{l|}{A91}            & \multicolumn{1}{l|}{B17}            & \multicolumn{1}{l|}{A262}           & \multicolumn{1}{l|}{B21}            & \multicolumn{1}{l|}{A257}           & B26            & A41            & B32            \\ \hline
\multicolumn{1}{|l|}{A141}           & \multicolumn{1}{l|}{B2}             & \multicolumn{1}{l|}{A31}            & \multicolumn{1}{l|}{B7}             & \multicolumn{1}{l|}{A74}            & \multicolumn{1}{l|}{B12}            & \multicolumn{1}{l|}{A161}           & \multicolumn{1}{l|}{B17}            & \multicolumn{1}{l|}{A272}           & \multicolumn{1}{l|}{B21}            & \multicolumn{1}{l|}{A270}           & B26            & A76            & B32            \\ \hline
\multicolumn{1}{|l|}{A176}           & \multicolumn{1}{l|}{B2}             & \multicolumn{1}{l|}{A66}            & \multicolumn{1}{l|}{B7}             & \multicolumn{1}{l|}{A96}            & \multicolumn{1}{l|}{B12}            & \multicolumn{1}{l|}{A196}           & \multicolumn{1}{l|}{B17}            & \multicolumn{1}{l|}{A297}           & \multicolumn{1}{l|}{B21}            & \multicolumn{1}{l|}{A327}           & B26            & A146           & B32            \\ \hline
\multicolumn{1}{|l|}{A211}           & \multicolumn{1}{l|}{B2}             & \multicolumn{1}{l|}{A101}           & \multicolumn{1}{l|}{B7}             & \multicolumn{1}{l|}{A131}           & \multicolumn{1}{l|}{B12}            & \multicolumn{1}{l|}{A231}           & \multicolumn{1}{l|}{B17}            & \multicolumn{1}{l|}{A332}           & \multicolumn{1}{l|}{B21}            & \multicolumn{1}{l|}{A46}            & B27            & A181           & B32            \\ \hline
\multicolumn{1}{|l|}{A246}           & \multicolumn{1}{l|}{B2}             & \multicolumn{1}{l|}{A111}           & \multicolumn{1}{l|}{B7}             & \multicolumn{1}{l|}{A166}           & \multicolumn{1}{l|}{B12}            & \multicolumn{1}{l|}{A266}           & \multicolumn{1}{l|}{B17}            & \multicolumn{1}{l|}{A16}            & \multicolumn{1}{l|}{B22}            & \multicolumn{1}{l|}{A81}            & B27            & A216           & B32            \\ \hline
\multicolumn{1}{|l|}{A256}           & \multicolumn{1}{l|}{B2}             & \multicolumn{1}{l|}{A136}           & \multicolumn{1}{l|}{B7}             & \multicolumn{1}{l|}{A201}           & \multicolumn{1}{l|}{B12}            & \multicolumn{1}{l|}{A301}           & \multicolumn{1}{l|}{B17}            & \multicolumn{1}{l|}{A51}            & \multicolumn{1}{l|}{B22}            & \multicolumn{1}{l|}{A116}           & B27            & A251           & B32            \\ \hline
\multicolumn{1}{|l|}{A281}           & \multicolumn{1}{l|}{B2}             & \multicolumn{1}{l|}{A171}           & \multicolumn{1}{l|}{B7}             & \multicolumn{1}{l|}{A236}           & \multicolumn{1}{l|}{B12}            & \multicolumn{1}{l|}{A336}           & \multicolumn{1}{l|}{B17}            & \multicolumn{1}{l|}{A86}            & \multicolumn{1}{l|}{B22}            & \multicolumn{1}{l|}{A151}           & B27            & A321           & B32            \\ \hline
\multicolumn{1}{|l|}{A316}           & \multicolumn{1}{l|}{B2}             & \multicolumn{1}{l|}{A206}           & \multicolumn{1}{l|}{B7}             & \multicolumn{1}{l|}{A271}           & \multicolumn{1}{l|}{B12}            & \multicolumn{1}{l|}{A20}            & \multicolumn{1}{l|}{B18}            & \multicolumn{1}{l|}{A121}           & \multicolumn{1}{l|}{B22}            & \multicolumn{1}{l|}{A186}           & B27            & A5             & B33            \\ \hline
\multicolumn{1}{|l|}{A35}            & \multicolumn{1}{l|}{B3}             & \multicolumn{1}{l|}{A241}           & \multicolumn{1}{l|}{B7}             & \multicolumn{1}{l|}{A306}           & \multicolumn{1}{l|}{B12}            & \multicolumn{1}{l|}{A55}            & \multicolumn{1}{l|}{B18}            & \multicolumn{1}{l|}{A156}           & \multicolumn{1}{l|}{B22}            & \multicolumn{1}{l|}{A221}           & B27            & A40            & B33            \\ \hline
\multicolumn{1}{|l|}{A70}            & \multicolumn{1}{l|}{B3}             & \multicolumn{1}{l|}{A276}           & \multicolumn{1}{l|}{B7}             & \multicolumn{1}{l|}{A341}           & \multicolumn{1}{l|}{B12}            & \multicolumn{1}{l|}{A90}            & \multicolumn{1}{l|}{B18}            & \multicolumn{1}{l|}{A191}           & \multicolumn{1}{l|}{B22}            & \multicolumn{1}{l|}{A291}           & B27            & A85            & B33            \\ \hline
\multicolumn{1}{|l|}{A105}           & \multicolumn{1}{l|}{B3}             & \multicolumn{1}{l|}{A286}           & \multicolumn{1}{l|}{B7}             & \multicolumn{1}{l|}{A351}           & \multicolumn{1}{l|}{B12}            & \multicolumn{1}{l|}{A125}           & \multicolumn{1}{l|}{B18}            & \multicolumn{1}{l|}{A226}           & \multicolumn{1}{l|}{B22}            & \multicolumn{1}{l|}{A326}           & B27            & A120           & B33            \\ \hline
\multicolumn{1}{|l|}{A140}           & \multicolumn{1}{l|}{B3}             & \multicolumn{1}{l|}{A311}           & \multicolumn{1}{l|}{B7}             & \multicolumn{1}{l|}{A25}            & \multicolumn{1}{l|}{B13}            & \multicolumn{1}{l|}{A135}           & \multicolumn{1}{l|}{B18}            & \multicolumn{1}{l|}{A261}           & \multicolumn{1}{l|}{B22}            & \multicolumn{1}{l|}{A10}            & B28            & A145           & B33            \\ \hline
\multicolumn{1}{|l|}{A175}           & \multicolumn{1}{l|}{B3}             & \multicolumn{1}{l|}{A346}           & \multicolumn{1}{l|}{B7}             & \multicolumn{1}{l|}{A60}            & \multicolumn{1}{l|}{B13}            & \multicolumn{1}{l|}{A160}           & \multicolumn{1}{l|}{B18}            & \multicolumn{1}{l|}{A296}           & \multicolumn{1}{l|}{B22}            & \multicolumn{1}{l|}{A45}            & B28            & A155           & B33            \\ \hline
\multicolumn{1}{|l|}{A210}           & \multicolumn{1}{l|}{B3}             & \multicolumn{1}{l|}{A30}            & \multicolumn{1}{l|}{B8}             & \multicolumn{1}{l|}{A95}            & \multicolumn{1}{l|}{B13}            & \multicolumn{1}{l|}{A170}           & \multicolumn{1}{l|}{B18}            & \multicolumn{1}{l|}{A331}           & \multicolumn{1}{l|}{B22}            & \multicolumn{1}{l|}{A80}            & B28            & A180           & B33            \\ \hline
\multicolumn{1}{|l|}{A280}           & \multicolumn{1}{l|}{B3}             & \multicolumn{1}{l|}{A65}            & \multicolumn{1}{l|}{B8}             & \multicolumn{1}{l|}{A130}           & \multicolumn{1}{l|}{B13}            & \multicolumn{1}{l|}{A195}           & \multicolumn{1}{l|}{B18}            & \multicolumn{1}{l|}{A15}            & \multicolumn{1}{l|}{B23}            & \multicolumn{1}{l|}{A115}           & B28            & A215           & B33            \\ \hline
\multicolumn{1}{|l|}{A315}           & \multicolumn{1}{l|}{B3}             & \multicolumn{1}{l|}{A100}           & \multicolumn{1}{l|}{B8}             & \multicolumn{1}{l|}{A165}           & \multicolumn{1}{l|}{B13}            & \multicolumn{1}{l|}{A230}           & \multicolumn{1}{l|}{B18}            & \multicolumn{1}{l|}{A50}            & \multicolumn{1}{l|}{B23}            & \multicolumn{1}{l|}{A150}           & B28            & A250           & B33            \\ \hline
\multicolumn{1}{|l|}{A350}           & \multicolumn{1}{l|}{B3}             & \multicolumn{1}{l|}{A205}           & \multicolumn{1}{l|}{B8}             & \multicolumn{1}{l|}{A200}           & \multicolumn{1}{l|}{B13}            & \multicolumn{1}{l|}{A240}           & \multicolumn{1}{l|}{B18}            & \multicolumn{1}{l|}{A62}            & \multicolumn{1}{l|}{B23}            & \multicolumn{1}{l|}{A185}           & B28            & A320           & B33            \\ \hline
\multicolumn{1}{|l|}{A34}            & \multicolumn{1}{l|}{B4}             & \multicolumn{1}{l|}{A275}           & \multicolumn{1}{l|}{B8}             & \multicolumn{1}{l|}{A245}           & \multicolumn{1}{l|}{B13}            & \multicolumn{1}{l|}{A265}           & \multicolumn{1}{l|}{B18}            & \multicolumn{1}{l|}{A190}           & \multicolumn{1}{l|}{B23}            & \multicolumn{1}{l|}{A220}           & B28            & A4             & B34            \\ \hline
\multicolumn{1}{|l|}{A69}            & \multicolumn{1}{l|}{B4}             & \multicolumn{1}{l|}{A285}           & \multicolumn{1}{l|}{B8}             & \multicolumn{1}{l|}{A305}           & \multicolumn{1}{l|}{B13}            & \multicolumn{1}{l|}{A300}           & \multicolumn{1}{l|}{B18}            & \multicolumn{1}{l|}{A225}           & \multicolumn{1}{l|}{B23}            & \multicolumn{1}{l|}{A255}           & B28            & A39            & B34            \\ \hline
\multicolumn{1}{|l|}{A104}           & \multicolumn{1}{l|}{B4}             & \multicolumn{1}{l|}{A310}           & \multicolumn{1}{l|}{B8}             & \multicolumn{1}{l|}{A340}           & \multicolumn{1}{l|}{B13}            & \multicolumn{1}{l|}{A335}           & \multicolumn{1}{l|}{B18}            & \multicolumn{1}{l|}{A235}           & \multicolumn{1}{l|}{B23}            & \multicolumn{1}{l|}{A290}           & B28            & A109           & B34            \\ \hline
\multicolumn{1}{|l|}{A139}           & \multicolumn{1}{l|}{B4}             & \multicolumn{1}{l|}{A345}           & \multicolumn{1}{l|}{B8}             & \multicolumn{1}{l|}{A37}            & \multicolumn{1}{l|}{B14}            & \multicolumn{1}{l|}{A19}            & \multicolumn{1}{l|}{B19}            & \multicolumn{1}{l|}{A260}           & \multicolumn{1}{l|}{B23}            & \multicolumn{1}{l|}{A325}           & B28            & A144           & B34            \\ \hline
\multicolumn{1}{|l|}{A174}           & \multicolumn{1}{l|}{B4}             & \multicolumn{1}{l|}{A29}            & \multicolumn{1}{l|}{B9}             & \multicolumn{1}{l|}{A59}            & \multicolumn{1}{l|}{B14}            & \multicolumn{1}{l|}{A89}            & \multicolumn{1}{l|}{B19}            & \multicolumn{1}{l|}{A295}           & \multicolumn{1}{l|}{B23}            & \multicolumn{1}{l|}{A9}             & B29            & A179           & B34            \\ \hline
\multicolumn{1}{|l|}{A209}           & \multicolumn{1}{l|}{B4}             & \multicolumn{1}{l|}{A64}            & \multicolumn{1}{l|}{B9}             & \multicolumn{1}{l|}{A94}            & \multicolumn{1}{l|}{B14}            & \multicolumn{1}{l|}{A124}           & \multicolumn{1}{l|}{B19}            & \multicolumn{1}{l|}{A330}           & \multicolumn{1}{l|}{B23}            & \multicolumn{1}{l|}{A44}            & B29            & A214           & B34            \\ \hline
\multicolumn{1}{|l|}{A222}           & \multicolumn{1}{l|}{B4}             & \multicolumn{1}{l|}{A99}            & \multicolumn{1}{l|}{B9}             & \multicolumn{1}{l|}{A129}           & \multicolumn{1}{l|}{B14}            & \multicolumn{1}{l|}{A159}           & \multicolumn{1}{l|}{B19}            & \multicolumn{1}{l|}{A14}            & \multicolumn{1}{l|}{B24}            & \multicolumn{1}{l|}{A54}            & B29            & A249           & B34            \\ \hline
\multicolumn{1}{|l|}{A244}           & \multicolumn{1}{l|}{B4}             & \multicolumn{1}{l|}{A112}           & \multicolumn{1}{l|}{B9}             & \multicolumn{1}{l|}{A164}           & \multicolumn{1}{l|}{B14}            & \multicolumn{1}{l|}{A194}           & \multicolumn{1}{l|}{B19}            & \multicolumn{1}{l|}{A24}            & \multicolumn{1}{l|}{B24}            & \multicolumn{1}{l|}{A79}            & B29            & A284           & B34            \\ \hline
\multicolumn{1}{|l|}{A254}           & \multicolumn{1}{l|}{B4}             & \multicolumn{1}{l|}{A134}           & \multicolumn{1}{l|}{B9}             & \multicolumn{1}{l|}{A199}           & \multicolumn{1}{l|}{B14}            & \multicolumn{1}{l|}{A229}           & \multicolumn{1}{l|}{B19}            & \multicolumn{1}{l|}{A49}            & \multicolumn{1}{l|}{B24}            & \multicolumn{1}{l|}{A114}           & B29            & A319           & B34            \\ \hline
\multicolumn{1}{|l|}{A279}           & \multicolumn{1}{l|}{B4}             & \multicolumn{1}{l|}{A169}           & \multicolumn{1}{l|}{B9}             & \multicolumn{1}{l|}{A234}           & \multicolumn{1}{l|}{B14}            & \multicolumn{1}{l|}{A239}           & \multicolumn{1}{l|}{B19}            & \multicolumn{1}{l|}{A84}            & \multicolumn{1}{l|}{B24}            & \multicolumn{1}{l|}{A149}           & B29            & A3             & B35            \\ \hline
\multicolumn{1}{|l|}{A314}           & \multicolumn{1}{l|}{B4}             & \multicolumn{1}{l|}{A204}           & \multicolumn{1}{l|}{B9}             & \multicolumn{1}{l|}{A269}           & \multicolumn{1}{l|}{B14}            & \multicolumn{1}{l|}{A264}           & \multicolumn{1}{l|}{B19}            & \multicolumn{1}{l|}{A119}           & \multicolumn{1}{l|}{B24}            & \multicolumn{1}{l|}{A184}           & B29            & A38            & B35            \\ \hline
\multicolumn{1}{|l|}{A349}           & \multicolumn{1}{l|}{B4}             & \multicolumn{1}{l|}{A274}           & \multicolumn{1}{l|}{B9}             & \multicolumn{1}{l|}{A339}           & \multicolumn{1}{l|}{B14}            & \multicolumn{1}{l|}{A299}           & \multicolumn{1}{l|}{B19}            & \multicolumn{1}{l|}{A154}           & \multicolumn{1}{l|}{B24}            & \multicolumn{1}{l|}{A219}           & B29            & A73            & B35            \\ \hline
\multicolumn{1}{|l|}{A33}            & \multicolumn{1}{l|}{B5}             & \multicolumn{1}{l|}{A309}           & \multicolumn{1}{l|}{B9}             & \multicolumn{1}{l|}{A23}            & \multicolumn{1}{l|}{B15}            & \multicolumn{1}{l|}{A334}           & \multicolumn{1}{l|}{B19}            & \multicolumn{1}{l|}{A189}           & \multicolumn{1}{l|}{B24}            & \multicolumn{1}{l|}{A289}           & B29            & A108           & B35            \\ \hline
\multicolumn{1}{|l|}{A68}            & \multicolumn{1}{l|}{B5}             & \multicolumn{1}{l|}{A28}            & \multicolumn{1}{l|}{B10}            & \multicolumn{1}{l|}{A58}            & \multicolumn{1}{l|}{B15}            & \multicolumn{1}{l|}{A344}           & \multicolumn{1}{l|}{B19}            & \multicolumn{1}{l|}{A224}           & \multicolumn{1}{l|}{B24}            & \multicolumn{1}{l|}{A324}           & B29            & A143           & B35            \\ \hline
\multicolumn{1}{|l|}{A103}           & \multicolumn{1}{l|}{B5}             & \multicolumn{1}{l|}{A63}            & \multicolumn{1}{l|}{B10}            & \multicolumn{1}{l|}{A128}           & \multicolumn{1}{l|}{B15}            & \multicolumn{1}{l|}{A18}            & \multicolumn{1}{l|}{B20}            & \multicolumn{1}{l|}{A259}           & \multicolumn{1}{l|}{B24}            & \multicolumn{1}{l|}{A8}             & B30            & A178           & B35            \\ \hline
\multicolumn{1}{|l|}{A113}           & \multicolumn{1}{l|}{B5}             & \multicolumn{1}{l|}{A98}            & \multicolumn{1}{l|}{B10}            & \multicolumn{1}{l|}{A163}           & \multicolumn{1}{l|}{B15}            & \multicolumn{1}{l|}{A53}            & \multicolumn{1}{l|}{B20}            & \multicolumn{1}{l|}{A294}           & \multicolumn{1}{l|}{B24}            & \multicolumn{1}{l|}{A43}            & B30            & A213           & B35            \\ \hline
\multicolumn{1}{|l|}{A138}           & \multicolumn{1}{l|}{B5}             & \multicolumn{1}{l|}{A133}           & \multicolumn{1}{l|}{B10}            & \multicolumn{1}{l|}{A198}           & \multicolumn{1}{l|}{B15}            & \multicolumn{1}{l|}{A88}            & \multicolumn{1}{l|}{B20}            & \multicolumn{1}{l|}{A304}           & \multicolumn{1}{l|}{B24}            & \multicolumn{1}{l|}{A78}            & B30            & A248           & B35            \\ \hline
\multicolumn{1}{|l|}{A148}           & \multicolumn{1}{l|}{B5}             & \multicolumn{1}{l|}{A168}           & \multicolumn{1}{l|}{B10}            & \multicolumn{1}{l|}{A233}           & \multicolumn{1}{l|}{B15}            & \multicolumn{1}{l|}{A123}           & \multicolumn{1}{l|}{B20}            & \multicolumn{1}{l|}{A329}           & \multicolumn{1}{l|}{B24}            & \multicolumn{1}{l|}{A126}           & B30            & A283           & B35            \\ \hline
                                     &                                     &                                     &                                     &                                     &                                     &                                     &                                     &                                     &                                     &                                     &                & A13            & B25            \\ \cline{13-14} 
\end{tabular}%
}
\end{table}


\begin{thebibliography}{}

\bibitem{RefJ}
J. Edmonds. "Paths, trees and 
flowers". Canadian Journal of Mathematics, 17: 449-467, 1965
\newline
\bibitem{RefJ}
S. Micali and V. V. Vazirani. "An $O(\sqrt{|V|} \cdot |E|)$ algorithm for finding maximum matching in general
graphs". In Proceedings of the twenty First annual IEEE
Symposium on Foundations of Computer Science,1980.
\newline
\bibitem{RefJ}
N. Blum. "A new approach to maximum matching in
general graphs". In Proc. 17th ICALP, volume 443 of
Lecture Notes in Computer Science, pages 586 - 597.
Springer-Verlag, 1990.
\newline
\bibitem{RefJ}
H. N. Gabow and R. E. Tarjan. "Faster scaling
algorithms for general graph matching problems". J.
ACM, 38(4):815 - 853, 1991.
\newline
\bibitem{RefJ}
Ryoji Kurata, Masahiro Goto, Atsushi
Iwasaki, and Makoto Yokoo. "Controlled school choice
with soft bounds and overlapping types". In AAAI Conference
on Artificial Intelligence (AAAI), 2015.
\newline
\bibitem{RefJ}
Dimitris Bertsimas, Vivek F Farias,
and Nikolaos Trichakis. "Fairness, efficiency, and flexibility
in organ allocation for kidney transplantation". Operations
Research, 61(1):73–87, 2013.
\newline
\bibitem{RefJ}
John Joseph Horton. "The effects of algorithmic
labor market recommendations: evidence from a field
experiment",To appear, Journal of Labor Economics, 2017.
\newline
\bibitem{RefJ}
E Krissinel and K Henrick.
"Secondary-structure matching (ssm), a new tool for fast
protein structure alignment in three dimensions". Acta
Crystallographica Section D: Biological Crystallography,
60(12):2256–2268, 2004.
\newline
\bibitem{RefJ}
Serge Belongie, Jitendra Malik, and
Jan Puzicha. "Shape matching and object recognition using
shape contexts". IEEE Transactions on Pattern Analysis
and Machine Intelligence, 24(4):509–522, 2002.
\newline
\bibitem{RefJ}
Liang Pang, Yanyan Lan, Jiafeng Guo,
Jun Xu, Shengxian Wan, and Xueqi Cheng. "Text matching
as image recognition". In AAAI Conference on Artificial
Intelligence (AAAI), 2016.
\newline
\bibitem{RefJ}
Gediminas Adomavicius
and YoungOk Kwon." Improving aggregate recommendation
diversity using ranking-based techniques". IEEE
Transactions on Knowledge and Data Engineering
(TKDE), 24(5):896–911, 2012.
\newline
\bibitem{RefJ}
Chaofeng Sha, Xiaowei Wu, and Junyu
Niu. "A framework for recommending relevant and diverse
items". In Proceedings of the International Joint Conference
on Artificial Intelligence (IJCAI), 2016.
\newline
\bibitem{RefJ}
Azin Ashkan, Branislav Kveton,
Shlomo Berkovsky, and Zheng Wen. "Optimal greedy
diversity for recommendation". In Proceedings of the
International Joint Conference on Artificial Intelligence
(IJCAI), pages 1742–1748, 2015.
\newline
\bibitem{RefJ}
A. Fischer and C. Igel, 
"An Introduction to Restricted Boltzmann Machines,
in Progress in Pattern Recognition, Image Analysis, 
Computer Vision, and Applications", vol. 7441 of 
Lecture Notes in Computer Science, pp. 14–36, 
Springer Berlin Heidelberg, Berlin, Heidelberg, 2012. 
\newline
\bibitem{RefJ}
Hinton, G. E. "A Practical Guide to Training Restricted
Boltzmann Machines". Technical Report, Department of Computer Science, University of Toronto,
2010.
\newline
\bibitem{RefJ}
H.W. Kuhn, "On the origin of the Hungarian method for the assignment problem, in J.K. Lenstra, A.H.G. Rinnooy Kan, A. Schrijver", 
History of Mathematical Programming, Amsterdam, North-Holland, 1991, pp. 77-81

\end{thebibliography}
\end{document}